
\input amstex 
\documentstyle{amsppt}

\magnification=\magstep1
\hoffset=0.3truecm  
\hfuzz=3pt
\vsize=21.5truecm

\font\tenbb=msbm10     
\font\sevenbb=msbm7
\newfam\bbfam \def\bb{\fam\bbfam\tenbb}  
\textfont\bbfam=\tenbb
\scriptfont\bbfam=\sevenbb

\chardef\at="40  

\def\A{{\roman A}}
\def\B{{\roman B}}
\def\C{{\roman C}}
\def\D{{\roman D}}
\def\natn{{\bb N}} 
\def\intz{{\bb Z}} 
\def\ratq{{\bb Q}} 
\def\smallintz{{\scriptstyle\bb Z}} 
\def\lcm{\hbox{{\rm lcm}}} 
\def\Int{\hbox{{\rm Int}}} 
\def\IntZ{\hbox{{\rm Int}}({\bb Z})} 
\def\IntZn{\hbox{{\rm Int}}({\bb Z}^n)} 

\def\x{\underline{x}}  
\def\xn{x_1,\ldots, x_n}  
  
\def\xr{x_1,\ldots, x_r}  
 
\def\({{\rm (}} \def\){{\rm )}} 

\line{\tenrm to appear in Comm.\ Algebra\hfill}
\vskip2.5truecm

\topmatter

\title 
Remarks on polynomial parametrization of sets of integer points
\endtitle

\author Sophie Frisch\endauthor

\address Institut f\"ur Mathematik C, Technische Universit\"at Graz,
{A-8010} Graz, Austria\endaddress

\email\nofrills frisch\@tugraz.at\endemail

\thanks
This note was written while the author was enjoying hospitality at
Universit\'e de Picardie, Amiens.
\endthanks



\abstract
If, for a subset $S$ of $\smallintz^k$, we compare the conditions
of being parametrizable $(a)$ by a single $k$-tuple of polynomials
with integer coefficients, $(b)$ by a single $k$-tuple of 
integer-valued polynomials and $(c)$ by finitely many $k$-tuples of
polynomials with integer coefficients (variables ranging through
the integers in each case), then $a\Rightarrow b$ (obviously),
$b\Rightarrow c$, and neither implication is reversible.
We give different characterizations of condition $(b)$.
Also, we show that every co-finite subset of $\smallintz^k$ is
parametrizable a single $k$-tuple of polynomials with integer
coefficients. 
\endabstract

\endtopmatter

\document

\noindent
If $f=(f_1,\ldots,f_k)\in(\intz[\xn])^k$ is a $k$-tuple of polynomials 
with integer coefficients in several variables, we call range or image of
$f$ the range of the function  $f\colon \intz^n\longrightarrow \intz^k$
defined by substitution of integers for the variables; and likewise for
a $k$-tuple of integer-valued polynomials 
$(f_1,\ldots,f_k)\in(\IntZn)^k$, where
$$\IntZn=\{g\in\ratq[\xn]\mid \forall a\in \intz^n: g(a)\in\intz\}.$$
If $S\subseteq\intz^k$ is the range of $f=(f_1,\ldots,f_k)$, we say that
$f$ parametrizes $S$.

We want to compare two kinds of polynomial parmetrization of sets of
integers or $k$-tuples of integers: by integer-valued polynomials and
by polynomials with integer coefficients.
Consider for instance the set of integer Pythagorean triples:
it takes two triples of polynomials with integer coefficients,
$\left(c(a^2-b^2),\; 2cab,\; c(a^2+b^2)\right)$ and
$\left(2cab,\; c(a^2-b^2),\; c(a^2+b^2)\right)$ to parametrize the
set of integer triples $(x,y,z)$ satisfying $x^2+y^2=z^2$, but
the same set can be parametrized by a single triple of integer-valued
polynomials \cite{2}.
Another reason for studying parametrization by
integer-valued polynomials are various sets of integers in number
theory and combinatorics that come parametrized by integer-valued
polynomials in a natural way, for example, the polygonal numbers 
$$p(n,k) = { {(n-2)k^2 - (n-4)k} \over 2}$$
where $p(n,k)$ represents the $k$-th $n$-gonal number \cite{3}. 

Now for our comparison of different kinds of polynomial
parametrization of sets of integer points.

\proclaim{Theorem}
For a set $S\subseteq \intz^k$ consider the conditions:
\smallskip

{\parskip=\smallskipamount
\item{\rm (A)}
$S$ is parmetrizable by a $k$-tuple of polynomials with
integer coefficients, i.e., there exists $f=(f_1,\ldots, f_k)$ in
$(\intz[\xn])^k$ \(for some $n$\) such that $S=f(\intz^n)$.
\item{\rm (B)}
$S$ is parmetrizable by a $k$-tuple of integer-valued polynomials,
i.e., there exists $g=(g_1,\ldots, g_k)$ in $(\Int(\intz^m))^k$
\(for some $m$\) such that $S=g(\intz^m)$.
\item{\rm (C)}
$S$ is a finite union of sets, each parametrizable by a $k$-tuple of 
polynomials with integer coefficients.
\item{\rm (D)}
$S$ is the set of integer $k$-tuples in the range of a 
$k$-tuple of polynomials with rational coefficients, as the
variables range through the integers, i.e., there exists
$h=(h_1,\ldots, h_k)$ in $(\ratq[\xr])^k$ \(for some $r$\) 
such that $S=h(\intz^r)\cap \intz^k$.
\smallskip
}

Then the following implications hold:

$$\matrix
A& & \cr
\Downarrow& & \cr
B&\Leftrightarrow &D\cr
\Downarrow& & \cr
C& & \cr
\endmatrix$$
\smallskip

and $\C \not\Rightarrow \B$, $\B\not \Rightarrow \A$.
\endproclaim
\medskip

Of the implications in the theorem,
$\A \Rightarrow \B$ and $\B\Rightarrow \D$ are trivial. 
We now show the nontrivial ones.

For $\D\Leftrightarrow \B$, we first construct,
for any $f\in \ratq[\xn]$, a parametrization of $f^{-1}(\intz)$ by
polynomials with integer coefficients, which we then plug into $f$ 
to obtain an integer-valued polynomial.

\proclaim{Lemma 1}
If $q_1,\ldots,q_r$ are powers of different primes and for each $i$,
$S_i$ is a union of residue classes of $q_i\intz^k$ in $\intz^k$ then 
$\bigcap_{i=1}^r S_i\subseteq \intz^k$ is parametrizable by a vector of 
polynomials with integer coefficients.
\endproclaim

\demo{Proof}
We will first parametrize a union of residue classes of $q\intz^k$ in
$\intz^k$ for a single prime power $q$. 
Let $a_0,\ldots,a_s\in \intz^k$ be representatives of the residue
classes in question, and let $t$ such that $2^t>s$. Expressing
$l\in\{0,1,\ldots,s\}$ in base $2$, we obtain a sequence of digits 
$[\kern1pt l\kern2pt]_2=(\varepsilon_0^{(l\/)},\ldots,\varepsilon_t^{(l\/)})$.
Let $m$ be a natural number such that $z^m$ is either congruent to $0$ or to
$1$ mod $q$ for every integer $z$.
Then
$$(qy_1,\ldots,qy_k) + \sum_{l=0}^s a_l\prod_{i=0}^t e_i^{(l\/)}(x_i),
\qquad\hbox{\rm with}\quad
e_i^{(l)}(x_i) = \left\{ 
\matrix x_i^m &\hbox{\rm\ if\ } \varepsilon_i^{(l\/)}=1\cr
1-x_i^m &\hbox{\rm\ if\ } \varepsilon_i^{(l\/)}=0\cr
\endmatrix
\right.
$$ 
parametrizes 
$\bigcup_{l=0}^s (q\intz^k + a_l)$.

Now let $q_1,\ldots,q_r$ be powers of different primes, and for 
$1\le i\le r$ let $S_i$ be a union of residue classes mod $q_i\intz^k$
parametrized by a polynomial vector $g_i$. By Chinese remainder theorem
there are $c_1,\ldots,c_r$ with $c_i\equiv 1$ mod $q_i$ and $c_i\equiv 0$
mod $q_j$ for $j\ne i$. We may choose $c_1,\ldots,c_r$ with
$\gcd(c_1,\ldots,c_r)=1$. (E.g.~by applying Dirichlet's theorem on primes
in arithmetic progressions to find primes $p_i\in b_i+q_i\intz$, where
$b_i$ is the inverse of $\prod_{j\ne i}q_j$ mod $q_i$, and setting 
$c_i=p_i\prod_{j\ne i}q_j$, with $p_1,\ldots,p_r$ different primes 
coprime to all $q_j$.)
Finally, we set $h=\sum_{i=1}^r c_i g_i$. 
Then $h$ parametrizes $\bigcap_{i=1}^r S_i$.
\qed \enddemo

\proclaim{Lemma 2 $(\B\Leftrightarrow D)$}
Let $S\subseteq \intz^k$. Then there exists a $k$-tuple of integer-valued 
polynomials whose range is $S$ if and only if there exists a $k$-tuple 
of polynomials with rational coefficients such that $S$ is the set of
integer points in its range \(as the variables range through the integers\).
\endproclaim

\demo{Proof}
The ``only if'' direction (that's $\B\Rightarrow \D$) is trivial. For
the other direction, $\D\Rightarrow \B$, first consider
the case $k=1$ of a single rational polynomial $f(\xn)=g(\xn)/c$
with $g(\xn)\in\intz[\xn]$ and $c\in\natn$.

Let $T=\{a\in \intz^n\mid f(a)\in\intz\}$. 
If $c=q_1\cdot\ldots\cdot q_r$ is the factorization of $c$ into
prime powers and $T_i=\{a\in \intz^n\mid g(a)\in q_i\intz\}$, then
$T=\bigcap_{i=1}^r T_i$.
For each $i$, $T_i$ is a union of residue classes of $q_i\intz^n$.
Hence $T$ is parametrizable by a polynomial
vector $(h_1,\ldots,h_n)\in\intz[\x]^n$.
Substituting $h_i$ for
$x_i$ in $f$, we obtain an integer-valued polynomial 
$p(\x)=f(h_1(\x),\ldots,h_n(\x))$ whose range is exactly 
the set of integers in the range of $f$.

In the case $k>1$, the argument for the set of integer points in the
range of a vector of rational polynomials $(f_1,\ldots,f_k)$, with
$f_i(\xn)=g_i(\xn)/c$, is similar, using 
$T_i=\{a\in \intz^n\mid \forall j:\; g_j(a)\in q_i\intz\}$.
\qed \enddemo

\proclaim{Lemma 3 $(\B\Rightarrow \C)$}
If a set $S\subseteq \intz^k$ is parametrizable by a single $k$-tuple
of integer-valued polynomials, it is parametrizable by a finite 
number of $k$-tuples of polynomials with integer coefficients.
\endproclaim

\demo{Proof}
First consider an integer-valued polynomial $f(x)$ in one variable of
degree $d$.
Recall that the binomial polynomials
${x\choose n}={{x(x-1)\ldots(x-n+1)}\over {n!}}$
form a basis of the $\intz$-module $\IntZ$, so that there exist
integers $a_0,\ldots,a_d$ with $f=\sum_{n=0}^d a_n {x\choose n}$.

It is easy to see that ${{cy+j}\choose n}\in \intz[y]$ for any $j$ whenever
$c$ is a common multiple of $1,2,\ldots,n$. Therefore for
$c=\lcm(1,2,\ldots,d)$ and arbitrary $j$,
$$f_j(y)=f(cy+j)=\sum_{n=0}^d a_n {{cy+j}\choose n}$$
is in $\intz[y]$; and clearly the image of $f$ is the union of
the images of $f_j$, for $j=0,\ldots,c-1$.

Regarding integer-valued poynomials in several variables, products of
binomial polynomials in one variable each
$\prod_{i=1}^n {x_i\choose n_i}$ form a basis of $\IntZn$
\cite{1, Prop.~XI.1.12}.
So, if $f\in \IntZn$ is of degree $d_i$ in $x_i$, and $c_i$ is a
common multiple of $1,2,\ldots, d_i$ then for each choice of
$j_1,\ldots, j_n$,
$f_{j_1,\ldots, j_n}=f(c_1y_1+j_1,\ldots, c_ny_n+j_n)$, as a 
$\intz$-linear combination of polynomials 
$\prod_{i=1}^n {c_iy_i+j_i\choose n_i}\in\intz[y_1,\ldots, y_n]$ is
a polynomial with integer coefficients and the image of $f$ is
the union of the images of the polynomials $f_{j_1,\ldots, j_n}$
with $0\le j_m< c_m$. 

The same argument shows that the image of a vector of polynomials
$(g_1,\ldots, g_k)$ in $(\IntZn)^k$ is the union of the images of
$c_1\cdot\ldots\cdot c_n$ vectors of polynomials in 
$(\intz[y_1,\ldots, y_n])^k$, where 
$c_i=\lcm(1,2,\ldots, d_i)$, $d_i$ denoting the highest degree 
of any $g_m$ in the $i$-th variable.
\qed \enddemo

\proclaim{Remark} 
$\B\not\Rightarrow A$ and $C\not\Rightarrow \B$:
Finite sets of more than one element witness
$C\not\Rightarrow \B$. 
The set of integer Pythagorean triples mentioned above
is parametrizable by a single triple of polynomials in $\Int(\intz^4)$,
but not by any triple of polynomials with integer coefficients in any
number of variables \cite{2} therefore $\B\not\Rightarrow A$.
\endproclaim

This completes the proof of the theorem.
The remainder of this note is devoted to the fact that every
co-finite set is parametrizable by a single vector of polynomials
with integer coefficients. (I was asked by L.~Vaserstein in 
connection with a remark in \cite{4} to publish a proof of this.)

\proclaim{Proposition}
Let $S\subseteq \intz^k$ such that $\intz^k\setminus S$ is finite.
Then there exists a $k$-tuple of polynomials with integer coefficients
whose range is $S$.
\endproclaim

\demo{Proof}
We may suppose that the complement of $S$ in $\intz^k$ is contained
in a cuboid 
$\prod_{i=1}^k [\kern1pt 0,n_i]=
[\kern1pt 0,n_1]\times\ldots\times[\kern1pt 0,n_k]$, 
with $n_i$ a non-negative integer for $1\le i\le k$.
We will first construct a polynomial vector whose image is
$\intz^k\setminus \prod_{i=1}^k [\kern1pt 0,n_i]$, by induction on $k$.

$k=1$:
for $n\ge 0$, the range of the polynomial $f$ below
is $\intz\setminus [\kern1pt 0,n]$:
$$f= - x_5^2(x_1^2+x_2^2+x_3^2+x_4^2 +1) +
(1- x_5^2)(x_1^2+x_2^2+x_3^2+x_4^2 +n+1).$$
Once we have a polynomial vector $(f_1,\ldots,f_{k-1})$ parametrizing
$\intz^{k-1}\setminus  \prod_{i=1}^{k-1} [\kern1pt 0,n_i]$ and a
polynomial $f$ with range $\intz\setminus [\kern1pt 0,n_k]$, we set
$$\displaylines{g_i= (1+{x_i}^2)(1-z^2)^{2m}f_i + z^2{x_i}\quad (1\le i< k)\cr
\hbox{\rm and}\quad g_k= (1+y^2)z^{2m}f + (1-z^2)y}$$
with $m$ sufficiently large, see below, and check that the range of
$(g_1,\ldots,g_{k})$ is
$\intz^k\setminus \prod_{i=1}^k [\kern1pt 0,n_i]$:
For $z=x_1=\ldots=x_{k-1}=0$ we get  $(f_1,\ldots,f_{k-1}, y)$, while
for $z\in \{1,-1\}$ and $y=0$, we have $(x_1,\ldots,x_{k-1}, f)$,
so that $(g_1,\ldots,g_{k})$ certainly covers the desired range.

Also, we stay within the desired range. Indeed, for $z=0$, the first
$k-1$ coordinates become $(1+{x_i}^2)f_i$, and their image lies
within the image of $(f_1,\ldots,f_{k-1})$, and for $z\in \{1,-1\}$ the
last coordinate is $(1+y^2)f$, whose image is contained in the image
of $f$.

Let $n=\max_i\{n_i\}$. 
By choosing $m$ sufficiently large such that 
$$|(1+x^2)(1-z^2)^{2m}|> |z^2x| +n
\quad \hbox{\rm and}\quad
|(1+y^2)z^{2m}|> |(1-z^2)y| +n$$
for all $z$ with $|z|\ge 2$ and all values of $x$ and $y$,
we make sure that $(g_1,\ldots,g_k)$ stays
within the desired range also for $|z|\ge 2$.

Having constructed a polynomial vector with range 
$\intz^k\setminus \prod_{i=1}^k [\kern1pt 0,n_i]$,
we can add additional values to the range, one by one, as follows.

If $g=(g_1,\ldots,g_k)$ is a polynomial vector whose image contains 
$\intz^k\setminus \prod_{i=1}^k [\kern1pt 0,n_i]$, but does not
contain $0\in\intz^k$, and $c$ is in $\prod_{i=1}^k [\kern1pt 0,n_i]$, let
$$h= w^{2t}g + (1-w^2)c,$$
with $t$ such that $2^{2t-2}>\max_i\{n_i\}$ then
the range of $h$ is exactly the range of $g$ together with the
(possibly additional) value $c$. If the value $c=0\in\intz^k$ is to be
added to the range of $g$, it must be added last.
\qed \enddemo
\goodbreak

\refstyle{C}

\enddocument